\documentclass[10pt,a4paper,twoside]{article}
\usepackage{amsmath}
\usepackage{amssymb}
\usepackage{amsthm}
\usepackage{geometry}
\usepackage{graphicx}
\usepackage{subfig}

\date{}
\title{\textbf{A Remark on the Stability of Peakons for the Degasperis-Procesi Equation}}
\author{\textbf{Andr\'e Kabakouala}\\
L.M.P.T., U.F.R Sciences et Techniques, Universit\'e de Tours, Parc Grandmont,\\
37200 Tours, France.\\
\\
Andre.Kabakouala@lmpt.univ-tours.fr}

\begin{document}

\theoremstyle{plain}
\newtheorem{Theo}{Theorem}[section]
\newtheorem{Pro}{Proposition}[section]
\newtheorem{Lem}{Lemma}[section]
\newtheorem{Cor}{Corollary}[section]
\newtheorem{Def}{Definition}[section]
\newtheorem{Pre}{Proof}
\newtheorem{Rem}{Remark}[section]
\newtheorem{Ex}{Example}[section]
\newtheorem*{Ack}{Acknowledgements}
\newtheorem{Cost}{Construction}[section]
\numberwithin{equation}{section}

\maketitle

\begin{abstract}
In this paper, we present a new argument (see Lemma \ref{Lemma PP2}) that allows us to simplify the proof of stability of peakons established in Lin and Liu (2009) (Theorem 1.1).
\end{abstract}
\vspace{1cm}

\section{Introduction}
In this paper, we consider the Degasperis-Procesi equation (DP)
\begin{equation}
u_{t}-u_{txx}+4uu_{x}=3u_{x}u_{xx}+uu_{xxx},
~~(t,x)\in\mathbb{R}_{+}\times\mathbb{R},
\label{1.1}
\end{equation}
with $u(0)=u_{0}\in L^{2}(\mathbb{R})$ and $(1-\partial^{2}_{x})u_{0}\in\mathcal{M}^{+}(\mathbb{R})$.

The DP equation is completely integrable (see \cite{MR2001531}) and has been proved to be physically relevant for water waves (see \cite{MR2481064}). It   possesses, among others,  the following conservation laws
\begin{equation}
E(u)=\int_{\mathbb{R}}yv=\int_{\mathbb{R}}\left(4v^{2}+5v^{2}_{x}+v^{2}_{xx}\right),~~
F(u)=\int_{\mathbb{R}}u^{3}=\int_{\mathbb{R}}\left(-v^{3}_{xx}+12vv^{2}_{xx}-48v^{2}v_{xx}+64v^{3}\right),
\label{i}
\end{equation}
where $y=(1-\partial^{2}_{x})u$ and $v=(4-\partial^{2}_{x})^{-1}u$. One can notice that the conservation law 
$E(\cdot)$ is equivalent to $\|\cdot\|^{2}_{L^{2}(\mathbb{R})}$. Indeed, using integration by parts (we assume that $u(\pm\infty)=v(\pm\infty)=v_{x}(\pm\infty)=0$), it holds
\begin{equation}
\|u\|^{2}_{L^{2}(\mathbb{R})}=\int_{\mathbb{R}}u^{2}=\int_{\mathbb{R}}(4v-v_{xx})^{2}
=\int_{\mathbb{R}}\left(16v^{2}+8v^{2}_{x}+v^{2}_{xx}\right)\sim E(u).
\label{1.5}
\end{equation}
In the sequel we will denote 
\begin{equation}
\|u\|_{\mathcal{H}}=\sqrt{E(u)}.
\label{1.5b}
\end{equation}

Applying $(1-\partial^{2}_{x})^{-1}(\cdot)$ to \eqref{1.1}, we obtain
\begin{equation}\label{int1}
u_{t}+\frac{1}{2}\partial_{x} u^{2}+\frac{3}{2}(1-\partial^{2}_{x})^{-1}\partial_{x}u^{2}=0,
~~(t,x)\in\mathbb{R}_{+}\times\mathbb{R}.
\end{equation}
In this form, the DP equation admits explicit solitary waves called \textit{peakons} (see \cite{MR2001531}) that are defined by
\begin{equation}
u(t,x)=\varphi_{c}(x-ct)=c\varphi(x-ct)
=ce^{-|x-ct|},
~~c\in\mathbb{R}^{*},
~~(t,x)\in\mathbb{R}_{+}\times\mathbb{R}.
\label{1.7}
\end{equation}

Our goal is to simplify the proof given in \cite{MR2460268} of the stability of a single peakon for the DP equation. Recall that the proof of the stability for the Camassa-Holm equation (CH) in \cite{MR1737505} follows from two integral relations between two conservation laws of CH, $\max_{\mathbb{R}} u$ and functions related to $u$. In \cite{MR2460268} the proof is more complicated, since all the local maxima and minima of $v=(4-\partial^{2}_{x})^{-1}u$ are involved in the relations. In this paper we present a simplification of this proof, where only the maximum of $v$ is involved in the relations. Our proof is thus closer to the proof for CH in \cite{MR1737505}. The main idea is the following: since $u$ is $L^{2}$-close to the peakon $\varphi_{c}(\cdot-\xi)$, for some $\xi\in\mathbb{R}$, and $(1-\partial^{2}_{x})u\in\mathcal{M}^{+}(\mathbb{R})$, it is easy to check that $u$ is actually $C^{0}$-close to the peakon, and thus $v$ is $C^{2}$-close to the \textit{smooth-peakon}:
\begin{equation}\label{smooth1}
\rho_{c}(x-\xi)=(4-\partial^{2}_{x})^{-1}\varphi_{c}(x-\xi)
=\frac{c}{3}e^{-|x-\xi|}-\frac{c}{6}e^{-2|x-\xi|},
~~x\in\mathbb{R}.
\end{equation}
First, since $\rho_{c}$, $\rho'_{c}$ and $\rho''_{c}$ are very small with respect to the amplitude $c$ outside of the interval 
$\Theta_0=[-6.7,6.7]$, we can restrict ourself to study $v$ on $\Theta_\xi = [\xi-6.7,\xi+6.7]$. Now we observe that $\rho''_{c}$ has strictly negative values in the interval $\mathcal{V}_0=[-\text{ln}\sqrt{2},\text{ln}\sqrt{2}]$, with $\rho_c'$ strictly positive on $[-6.7,-\text{ln}\sqrt{2}]$ and $\rho_c' $ strictly negative on 
$[\text{ln}\sqrt{2},6.7]$. This forces $v_{x}$ to change sign only one time on $\Theta_\xi $, and thus $v$ has only one local extremum (which is a maximum) on $\Theta_\xi$. This fact will considerably simplify the proof of the stability.

\section{Preliminaries}\label{Section 2}
In this section, we briefly recall the global well-posedness results for the Cauchy problem of the DP equation 
(see \cite{MR2271927} and \cite{MR2249792}), and its consequences. For $I$ a finite or infinite time interval of $\mathbb{R}_{+}$, we denote by $\mathcal{X}(I)$ the function space \footnote{$ W^{1,1}(\mathbb R) $ is the space of $ L^1(\mathbb R) $ functions with derivatives in $ L^1(\mathbb R) $ and $ BV(\mathbb R) $ is the space of function with bounded variation.}
\begin{equation}
\mathcal{X}(I)=\left\lbrace u\in C\left(I;H^{1}(\mathbb{R})\right)\cap L^{\infty}\left(I;W^{1,1}(\mathbb{R})\right),~ u_{x}\in L^{\infty}\left(I;BV(\mathbb{R})\right)\right\rbrace.
\label{1.8}
\end{equation}

\begin{Theo}[Global Weak Solution; See \cite{MR2271927} and \cite{MR2249792}]\label{Theoreme 2.1}
Assume that  $u_{0}\in L^{2}(\mathbb{R})$ with  $y_{0}=(1-\partial^{2}_{x})u_{0}\in\mathcal{M}^{+}(\mathbb{R})$. Then the DP equation has a unique global weak solution $u\in\mathcal{X}(\mathbb{R}_{+})$ such that 
\begin{equation}
y(t,\cdot)=(1-\partial^{2}_{x})u(t,\cdot)\in \mathcal{M}^{+}(\mathbb{R}),~~\forall t\in\mathbb{R}_+
\label{2.1}
\end{equation}
and
\begin{equation}
|u_{x}(t,x)|\le u(t,x),~~\forall(t,x)\in\mathbb{R}_{+}\times\mathbb{R}.
\label{2.1a}
\end{equation}
Moreover $ E(\cdot) $ and $ F(\cdot) $ are conserved by the flow.
\end{Theo}

\begin{Rem}[Control of $L^{\infty}$ Norm by $L^{2}$ Norm]\label{Remark 2.1}
\normalfont
 \eqref{2.1a} and the well-known Sobolev embedding of $ H^1(\mathbb{R}) $ into $ L^\infty(\mathbb{R}) $ lead to 
\begin{equation}
\|u\|_{L^{\infty}(\mathbb{R})}\le \frac{1}{\sqrt{2}} \|u\|_{H^{1}(\mathbb{R})}\le 
\|u\|_{L^{2}(\mathbb{R})}.
\label{2.3}
\end{equation}
\end{Rem}

\section{Stability of peakons}\label{Section 3a}
In this section, we present our simplification of the proof of stability of peakons for the DP equation.

\begin{Theo}[Stability of Peakons]\label{Theorem P}
Let $u\in \mathcal{X}([0,T[)$, with $0<T\le+\infty$, be a solution of the DP equation and 
$\varphi_{c}$ be the peakon defined in \eqref{1.7}, traveling to the right at the speed $c>0$. There exist $C>0$ and $\varepsilon_{0}>0$ only depending on the speed $c$, such that if
\begin{equation}
y_{0}=(1-\partial^{2}_{x})u_{0}\in\mathcal{M}^{+}(\mathbb{R})
\label{P1}
\end{equation}
and
\begin{equation}
\|u_{0}-\varphi_{c}\|_{\mathcal{H}}\le\varepsilon^{2},~~\text{with}~~0<\varepsilon<\varepsilon_{0},
\label{P2}
\end{equation}
then 
\begin{equation}
\|u(t,\cdot)-\varphi_{c}(\cdot-\xi(t))\|_{\mathcal{H}}\le C\sqrt{\varepsilon},~~\forall t\in[0,T[,
\label{P3}
\end{equation}
where $\xi(t)\in\mathbb{R}$ is the only point where the function 
$v(t,\cdot)=(4-\partial^{2}_{x})^{-1}u(t,\cdot)$ attains its maximum.
\end{Theo}

We first  recall that $E(u)\sim E(\varphi_{c})$ and $F(u)\sim F(\varphi_{c})$ in $\mathbb{R}$, if 
$u\sim\varphi_{c}$ in $L^{2}(\mathbb{R})$, with $y\in\mathcal{M}^{+}(\mathbb{R})$
(see for instance \cite{MR2460268} or \cite{AK}).

\begin{Lem}[Control of Distances Between Energies; See \cite{AK}]\label{Lemma P1}
Let $u\in L^{2}(\mathbb{R})$ with $y=(1-\partial^{2}_{x})u\in\mathcal{M}^{+}(\mathbb{R})$. 
If $\|u-\varphi_{c}\|_{\mathcal{H}}\le\varepsilon^2 $, then
\begin{equation}
|E(u)-E(\varphi_{c})|\le O(\varepsilon^2)
\label{P4}
\end{equation}
and 
\begin{equation}
|F(u)-F(\varphi_{c})|\le O(\varepsilon^2),
\label{P5}
\end{equation}
where $O(\cdot)$ only depends on the speed $c$.
\end{Lem}
To prove Theorem \ref{Theorem P}, by the conservation of $E(\cdot)$, $F(\cdot)$ and the continuity of the map $t\mapsto u(t)$ from $[0,T[$ to $\mathcal{H}$ (since $\mathcal{H}\simeq L^{2}$), it suffices to prove that for any function $u\in L^{2}(\mathbb{R})$ satisfying $y=(1-\partial^{2}_{x})u\in\mathcal{M}^{+}(\mathbb{R})$, \eqref{P4} and \eqref{P5}, if 
\begin{equation}
\inf_{z\in\mathbb{R}}\|u-\varphi_{c}(\cdot-z)\|_{\mathcal{H}}\le\varepsilon^{1/4},
\label{P6}
\end{equation}
then 
\begin{equation}
\|u-\varphi_{c}(\cdot-\xi)\|_{\mathcal{H}}\le C\sqrt{\varepsilon},
\label{P7}
\end{equation}
where $\xi\in\mathbb{R}$ is the only point of maximum of $v$.

\begin{figure}[ht]
\centering
\subfloat[$\rho(x)=\frac{1}{3}e^{-|x|}-\frac{1}{6}e^{-2|x|}$ profile.]
{\includegraphics[width=8cm, height=6cm]{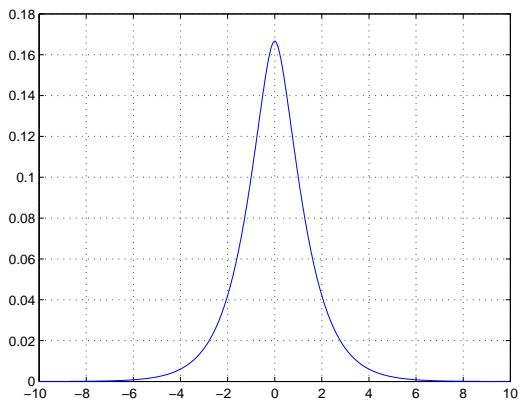}\label{F1}} \\
\subfloat[$\rho'(x)=\left(\frac{1}{3}e^{-|x|}-\frac{1}{3}e^{-2|x|}\right)_{x<0}
+\left(\frac{1}{3}e^{-2|x|}-\frac{1}{3}e^{-|x|}\right)_{x>0}$ profile.]
{\includegraphics[width=8cm, height=6cm]{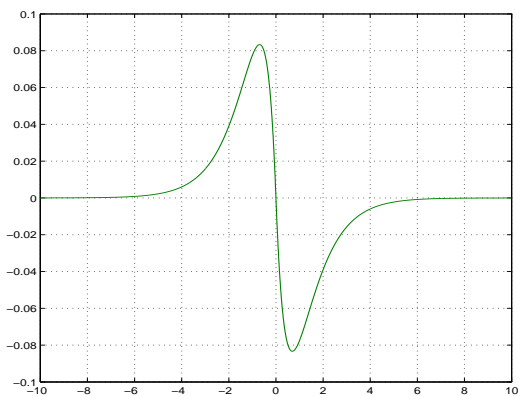}\label{F2}}
\subfloat[ $\rho''(x)=\frac{1}{3}e^{-|x|}-\frac{2}{3}e^{-2|x|}$ profile]
{\includegraphics[width=8cm, height=6cm]{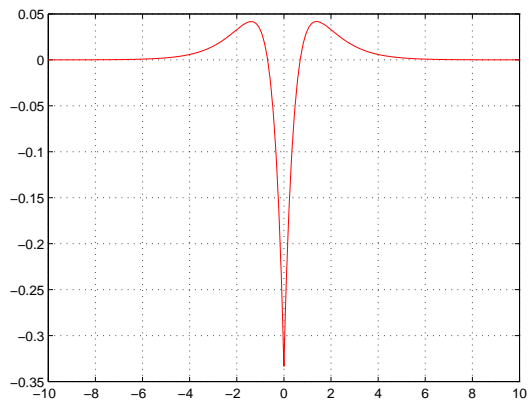}\label{F3}}
\caption{Variation of the smooth-peakon with the amplitude $1/6$ at initial time.}
\end{figure}

Let us present some important properties of smooth-peakons, defined in \eqref{smooth1}, which will play a crucial role in the proof of Theorem \ref{Theorem P}. The smooth-peakon $\rho_{c}$ belongs to $H^{3}(\mathbb{R})\hookrightarrow C^{2}(\mathbb{R})$ (by the Sobolev embedding) since $\varphi_{c}$ belongs to $H^{1}(\mathbb{R})$ (defined in \eqref{1.7}). It is a positive even function, which admits a single maximum $c/6$ at point $0$, and decays at infinity to $0$ (see Fig. \ref{F1}). Its derivative $\rho'_{c}$ belongs to $H^{2}(\mathbb{R})\hookrightarrow C^{1}(\mathbb{R})$, it is an odd function, which vanishes only at the origin and has  negative values on $[0,+\infty[$. It admits a  a single minimum $-c/12$ at point $\text{ln}2$ and tends at infinity to $0$ (see Fig. \ref{F2}). Its second derivative $\rho''_{c}$ belongs to $H^{1}(\mathbb{R})\hookrightarrow 
C^{0}(\mathbb{R})$, it is an even function, which vanishes at $\pm\text{ln}2$, takes positive values on 
$]-\infty,-\text{ln}2[\cup]\text{ln}2,+\infty[$ and negative values on $[-\text{ln}2,\text{ln}2]$. It admits a single minimum $-c/3$ at point $0$ and two maxima $c/24$ at points $\pm\text{ln}4$, and decays at infinity to $0$ (see Fig. \ref{F3}).

Next, we will need the following estimates.

\begin{Lem}[$C^{0}$, $C^{1}$ and $ C^2 $ Approximations]\label{Lemma PP1}
Let $u\in L^{2}(\mathbb{R})$ with $y=(1-\partial^{2}_{x})u\in\mathcal{M}^{+}(\mathbb{R})$. If 
$\|u-\varphi_{c}\|_{\mathcal{H}}\le \varepsilon^{1/4}$, then
\begin{equation}
\|u-\varphi_{c}\|_{C^{0}(\mathbb{R})}+ \|v-\rho_{c}\|_{C^{2}(\mathbb{R})}\le O(\varepsilon^{1/8})
\label{PP2}
\end{equation}
and 
\begin{equation}
\|v-\rho_{c}\|_{C^{1}(\mathbb{R})}\le O(\varepsilon^{1/4}).
\label{PP3}
\end{equation}
\end{Lem}
\noindent
\textbf{Proof.}
Let us begin with the second estimate. From the definition of $E(\cdot)$ and $\mathcal{H}$ 
(see respectively \eqref{i} and \eqref{1.5b}), one can see that $\|u\|_{\mathcal{H}}$ is equivalent to $\|v\|_{H^{2}(\mathbb{R})}$, since $\|v\|_{H^{2}(\mathbb{R})}\le\|u\|_{\mathcal{H}}\le 
5\|v\|_{H^{2}(\mathbb{R})}$. Then, assumption $u$ is $\mathcal{H}$-close to $\varphi_{c}$ implies that $v$ is $H^{2}$-close to $\rho_{c}$. Now, using the Sobolev embedding of $H^{2}(\mathbb{R})$ into $C^{1}(\mathbb{R})$, we deduce \eqref{PP3}.

For the first estimate,  note that the assumption $y=(1-\partial^{2}_{x})u\ge 0$ implies that 
$u=(1-\partial^{2}_{x})^{-1}y\ge 0$ and satisfies $|u_{x}|\le u$ on $\mathbb{R}$ (see \eqref{2.1a}). Then, applying triangular inequality, and using that $|\varphi'_{c}|=\varphi_{c}$ on $\mathbb{R}$ and \eqref{2.3}, we have 
\begin{align*}
\|u-\varphi_{c}\|_{H^{1}(\mathbb{R})}&\le\|u\|_{H^{1}(\mathbb{R})}
+\|\varphi_{c}\|_{H^{1}(\mathbb{R})}\\
&\le 2\|u\|_{L^{2}(\mathbb{R})}+2\|\varphi_{c}\|_{L^{2}(\mathbb{R})}\\
&\le 2\|u-\varphi_{c}\|_{L^{2}(\mathbb{R})}+4\|\varphi_{c}\|_{L^{2}(\mathbb{R})}\\
&\le O(\varepsilon^{1/4})+O(1),
\end{align*}
where $\|\varphi_{c}\|_{L^{2}(\mathbb{R})}=c$. Therefore, applying the Gagliardo-Nirenberg inequality and using that $\|u-\varphi_{c}\|_{\mathcal{H}}\le \varepsilon^{1/4}$ (with $\mathcal{H}\simeq L^{2}$), we obtain
\begin{align*}
\|u-\varphi_{c}\|_{C^{0}(\mathbb{R})}&\le 
\|u-\varphi_{c}\|^{1/2}_{L^{2}(\mathbb{R})}
\|u-\varphi_{c}\|^{1/2}_{H^{1}(\mathbb{R})}\\
&\le O(\varepsilon^{1/8})\left(O(\varepsilon^{1/8})+O(1)\right)\\
&\le O(\varepsilon^{1/8})\; .
\end{align*}
Finally to estimate the second term of the left-hand side of (\ref{PP2}), we first notice that  the continuity of $ (4-\partial_x^2) ^{-1}(\cdot) $ from $ H^s(\mathbb{R}) $ to $ H^{s+2}(\mathbb{R}) $ and the above estimates ensure that $\|v-\rho_c \|_{H^3} =O(1) $ and $\| v-\rho_c\|_{H^2} =O(\varepsilon^{1/4}) $ . These last estimates combined with  the Gagliardo-Nirenberg inequality yield the result as above.
\hfill $ \square $ \vspace*{2mm} 

The following lemma specifies the distance to minimize for stability.

\begin{Lem}[Quadratic Identity; See \cite{MR2460268}]\label{Lemma P2}
For any $u\in L^{2}(\mathbb{R})$ and $\xi\in\mathbb{R}$, it holds
\begin{equation}
E(u)-E(\varphi_{c})=\|u-\varphi_{c}(\cdot-\xi)\|^{2}_{\mathcal{H}}+4c\left(v(\xi)-\frac{c}{6}\right),
\label{P8}
\end{equation}
where $v=(4-\partial^{2}_{x})^{-1}u$ and $\rho_{c}(0)=c/6$.
\end{Lem}
\noindent
\textbf{Proof.} We follow the idea of Constantin and Strauss with the CH equation (see \cite{MR1737505}, Lemma $1$). We compute 
\begin{align}\label{P9}
E(u-\varphi_{c}(\cdot-\xi))&=E(u)+E(\varphi_{c})-2\left\langle(1-\partial^{2}_{x})\varphi_{c}(\cdot-\xi),
(4-\partial^{2}_{x})^{-1}u\right\rangle_{H^{-1},H^{1}}\nonumber\\
&=E(u)+E(\varphi_{c})-2\left\langle(1-\partial^{2}_{x})\varphi_{c}(\cdot-\xi),
v\right\rangle_{H^{-1},H^{1}},
\end{align}
where $\langle\cdot,\cdot\rangle_{H^{-1},H^{1}}$ denotes the duality bracket $H^{-1}(\mathbb{R})$, $H^{1}(\mathbb{R})$.
Now, using the definition of $\varphi'_{c}(\cdot-\xi)$ and integration by parts, we have  
\begin{align}\label{P10}
\left\langle(1-\partial^{2}_{x})\varphi_{c}(\cdot-\xi),v\right\rangle_{H^{-1},H^{1}}&=
\int_{\mathbb{R}}v\varphi_{c}(\cdot-\xi)+
\int_{\mathbb{R}}v_{x}\varphi'_{c}(\cdot-\xi)\nonumber\\
&=\int_{\mathbb{R}}v\varphi_{c}(\cdot-\xi)+
\int_{-\infty}^{\xi}v_{x}\varphi_{c}(\cdot-\xi)
-\int_{\xi}^{+\infty}v_{x}\varphi_{c}(\cdot-\xi)\nonumber\\
&=2cv(\xi).
\end{align}
Recalling  that the energy of peakons is given by 
\begin{align}\label{P11}
E(\varphi_{c})&=
\left\langle(1-\partial^{2}_{x})\varphi_{c},(4-\partial^{2}_{x})^{-1}\varphi_{c}
\right\rangle_{H^{-1},H^{1}}
=\int_{\mathbb{R}}\rho_{c}\varphi_{c}+\int_{\mathbb{R}}\rho'_{c}\varphi'_{c}\nonumber\\
&=\int_{\mathbb{R}}\rho_{c}\varphi_{c}
+\int_{-\infty}^{0}\rho'_{c}\varphi_{c}
-\int_{0}^{+\infty}\rho'_{c}\varphi_{c}
=2c\rho_{c}(0)
=\frac{c^{2}}{3} \; , 
\end{align}
we  obtain the lemma. \hfill $\square$ \vspace*{2mm}

Now we will study carefully the local extrema of $v=(4-\partial^{2}_{x})^{-1}u$. Let $u\in L^{2}(\mathbb{R})$ with $y = (1-\partial^{2}_{x})u\in\mathcal{M}^{+}(\mathbb{R})$, and assume that \eqref{P6} holds for some $z\in\mathbb{R}$. We consider the interval in which the mass of smooth-peakons is concentrated, and the interval in which the mass of second derivative of smooth-peakons is strictly negative. In the sequel of this paper, the notation $\alpha\simeq\beta$ means that $0.9\times\beta\le\alpha\le 1.1\times\beta$. We set, for any $ z\in \mathbb R $, 
\begin{equation}
\Theta_z=[z-6.7,z+6.7],~~\text{where}~~6.7\simeq\text{ln}\left(\frac{20}{20-\sqrt{399}}\right),
\label{P12}
\end{equation}
and
\begin{equation}\label{P13}
\mathcal{V}_z=\left[z-\text{ln}\sqrt{2},z+\text{ln}\sqrt{2}\right].
\end{equation}
One can clearly see that $\mathcal{V}_0$ is a subset of $\Theta_0$ (since $20/(20-\sqrt{399})>\sqrt{2}$). We chose the values $\pm 6.7$ such that $\rho_{c}(\pm 6.7)\simeq c/2400\simeq 4.1\times 
10^{-4}c$ as in \cite{AK}. Also, we have $\rho'_{c}(-6.7)=-\rho'_{c}(6.7)\simeq 4.1\times 10^{-4}c$ and 
$\rho''_{c}(\pm 6.7)\simeq 4.1\times 10^{-4}c$. Then $\rho_{c}$, 
$\rho'_{c}$ and $\rho''_{c}$ are very small with respect to the amplitude $c$ on $\mathbb{R}\setminus\Theta_0$. 

We claim the following result.

\begin{Lem}[Uniqueness of the Local Maximum]\label{Lemma PP2}
Let $u\in L^{2}(\mathbb{R})$ with $y=(1-\partial^{2}_{x})u\in\mathcal{M}^{+}(\mathbb{R})$, that satisfies 
\eqref{P6} for some $z\in\mathbb{R}$. There exists $\varepsilon_{0}>0$ only depending on the speed $c$, such that if $0<\varepsilon<\varepsilon_{0}$, then the function $v=(4-\partial^{2}_{x})^{-1}u$ admits a unique local extremum on $\Theta_z $. This extremum is a maximum, and it holds
\begin{equation}
v(x)\le\frac{c}{300},~~\forall x\in\mathbb{R}\setminus\Theta_z,
\label{P14}
\end{equation}
\begin{equation}
u(x)\le\frac{c}{300},~~\forall x\in\mathbb{R}\setminus\Theta_z.
\label{P15}
\end{equation}
\end{Lem}

\noindent
\textbf{Proof.} The key is to study the impact of the assumption $y\in\mathcal{M}^{+}(\mathbb{R})$ on $v$. First, let us show that $|v_{x}|\le 2v$ on $\mathbb{R}$. We recall that from the assumption $y\ge 0$, we have $u\ge 0$ and $v\ge 0$ on $\mathbb{R}$. According to the definition of $v$, we have for all $x\in\mathbb{R}$,
$$v(x)=\frac{e^{-2x}}{4}\int_{-\infty}^{x}e^{2x'}u(x')dx'
+\frac{e^{2x}}{4}\int_{x}^{+\infty}e^{-2x'}u(x')dx'$$
and 
$$v_{x}(x)=-\frac{e^{-2x}}{2}\int_{-\infty}^{x}e^{2x'}u(x')dx'
+\frac{e^{2x}}{2}\int_{x}^{+\infty}e^{-2x'}u(x')dx',$$
which yields
\begin{equation}
|v_{x}(x)|\le 2v(x),~~\forall x\in\mathbb{R}.
\label{P17}
\end{equation}

Second, let us show that $u\le 6v$ on $\mathbb{R}$.  Using the Fourier transform, one can check that
\begin{align} \label{P18}
(1-\partial^{2}_{x})^{-1}(4-\partial^{2}_{x})^{-1}(\cdot)
&=\mathcal{F}^{-1}\left[\frac{1}{3(1+\omega^{2})}-\frac{1}{3(4+\omega^{2})}\right](\cdot)\nonumber\\
&=\frac{1}{3}(1-\partial^{2}_{x})^{-1}(\cdot)
-\frac{1}{3}(4-\partial^{2}_{x})^{-1}(\cdot),
\end{align}
and one can rewrite $v$ as 
\begin{equation}\label{P19}
v=(4-\partial^{2}_{x})^{-1}(1-\partial^{2}_{x})^{-1}y=\frac{1}{3}(1-\partial^{2}_{x})^{-1}y
-\frac{1}{3}(4-\partial^{2}_{x})^{-1}y.
\end{equation}
Then for all $x\in\mathbb{R}$, 
\begin{align}\label{P20}
u(x)-6v(x)&=-(1-\partial^{2}_{x})^{-1}y(x)+2(4-\partial^{2}_{x})^{-1}y(x)\nonumber\\
&=-\frac{1}{2}\int_{\mathbb{R}}e^{-|x-x'|}y(x')dx'
+\frac{1}{2}\int_{\mathbb{R}}e^{-2|x-x'|}y(x')dx'\nonumber\\
&\le 0,
\end{align}
since $e^{-2|\cdot|}\le e^{-|\cdot|}$ on $\mathbb{R}$.

We are now ready to prove the uniqueness of local maxima in $ \Theta_z $.  Let us first study the sign of $v_{xx}$ on $\mathcal{V}_z$. One can easy check that for all $x\in\mathcal{V}_0$,
\begin{equation}\label{Valeur2}
\rho''_{c}(x)\le\frac{\sqrt{2}-2}{6}c.
\end{equation}
Then, combining \eqref{PP2} and \eqref{Valeur2}, taking $0<\varepsilon<\varepsilon_{0}$ with $\varepsilon_{0}\ll 1$, we have for all $x\in\mathcal{V}_z$,
\begin{equation*}
v_{xx}(x)\le\frac{\sqrt{2}-2}{6}c+O(\varepsilon^{1/4})
\le\frac{\sqrt{2}-2}{600}c<0,
\end{equation*}
which implies that $v_{x}$ is strictly decreasing on $\mathcal{V}_z$. Let us study the sign of $v_{x}$ on $\Theta_z\setminus\mathcal{V}_z$. One can easily check that 
\begin{equation}\label{Valeur1}
\rho'_{c}\left(-\text{ln}\sqrt{2}\right)
=\frac{\sqrt{2}-1}{6}c~~\text{and}~~
\rho'_{c}\left(\text{ln}\sqrt{2}\right)
=-\frac{\sqrt{2}-1}{6}c,
\end{equation}
and that $\rho'_{c}(x)\ge  10^{-4}c$ for all $x\in[-6.7,-\text{ln}\sqrt{2}]$. Then using \eqref{PP3} and taking $0<\varepsilon<\varepsilon_{0}$ with $\varepsilon_{0}\ll 1$, we have $v_{x}(x)\ge 4\times 10^{-5}c>0$ for all $x\in[z-6.7,z-\text{ln}\sqrt{2}]$. Proceeding in the same way, we obtain $v_{x}(x)\le -4\times 10^{-5}c<0$ for all $x\in[z+\text{ln}\sqrt{2},z+6.7]$. Since $v_{x}$ is strictly decreasing on $\mathcal{V}_z$ and changes sign,  $v_{x}$ vanishes once on $\mathcal{V}_z$ and thus on  $\Theta_z $. Hence,  $v$ admits a single local extremum on $\Theta_z$, which is a maximum since $v_{xx}<0$ on $\mathcal{V}_z$.

Now, using that $\rho_{c}$ is increasing on $\mathbb{R}^-$, \eqref{PP3} and taking $0<\varepsilon<\varepsilon_{0}$ with $\varepsilon_{0}\ll 1$, it holds for all $x\in]-\infty,z-6.7[$,
$$v(x)=\rho_{c}(x-z)+O(\varepsilon^{1/4})
\le\frac{c}{2400}+O(\varepsilon^{1/4})
\le\frac{c}{300}.$$
Proceeding in the same way for $x\in]z+6.7,+\infty[$, we obtain \eqref{P14}.

Combining \eqref{PP2}, \eqref{P20} and proceeding as for the estimate \eqref{P14}, we get 
\eqref{P15}. Note that $\varphi_{c}(\pm 6.7)\simeq 1.2 \times 10^{-3}c$. This completes the proof of the lemma. \hfill $ \square $ \vspace*{2mm}

Under the assumptions of Lemma \ref{Lemma PP2},  $v$ has got a unique point of global maximum on $\mathbb{R}$. In the sequel of this section, we will denote by $\xi$ this point of global maximum and we set
$M=v(\xi)=\max_{x\in\mathbb{R}}v(x)$. The next two lemmas can be directly deduced from the similar lemmas established in \cite{MR2460268} (see also  \cite{AK}).

\begin{Lem}[Connection Between $E(\cdot)$ and $M^{2}$; See \cite{MR2460268}]\label{Lemma P4}
Let $u\in L^{2}(\mathbb{R})$ and $v=(4-\partial^{2}_{x})^{-1}u\in H^{2}(\mathbb{R})$. Define the function $g$ by
\begin{equation}
  g(x)=\left\{
    \begin{aligned}
     &2v+v_{xx}-3v_{x},~~x<\xi,\\
     &2v+v_{xx}+3v_{x},~~x>\xi.\\
    \end{aligned}
  \right.
  \label{GG1}
\end{equation}
Then it holds
\begin{equation}
\int_{\mathbb{R}}g^{2}(x)dx=E(u)-12M^{2}.
\label{GG2}
\end{equation}
\end{Lem}

\begin{Lem}[Connection Between $F(\cdot)$ and $M^{3}$; See \cite{MR2460268}]\label{Lemma P5}
Let $u\in L^{2}(\mathbb{R})$ and $v=(4-\partial^{2}_{x})^{-1}u\in H^{2}(\mathbb{R})$. Define the function $h$ by
\begin{equation}
  h(x)=\left\{
    \begin{aligned}
     &-v_{xx}-6v_{x}+16v,~~x<\xi,\\
     &-v_{xx}+6v_{x}+16v,~~x>\xi.\\
    \end{aligned}
  \right.
  \label{HH1}
\end{equation}
Then it holds
\begin{equation}
\int_{\mathbb{R}}h(x)g^{2}(x)dx=F(u)-144M^{3}.
\label{HH2}
\end{equation}
\end{Lem}

\textbf{Sketch of proof.} The proof of Lemmas \ref{Lemma P4}-\ref{Lemma P5} follows by direct computation, using integration by parts, with
$v_{x}(\xi)=0$ and $v(\pm\infty)=v_{x}(\pm\infty)=v_{xx}(\pm\infty)=0$. See \cite{MR2460268} (also \cite{AK})
to undersand the technique. 
\hfill $ \square $ \vspace*{2mm}

We can now connect the conservation laws.

\begin{Lem}[Connection Between $E(\cdot)$ and $F(\cdot)$]\label{Lemma P6}
Let $u\in L^{2}(\mathbb{R})$, with  $y=(1-\partial^{2}_{x})u \in {\mathcal M}^+({\mathbb R})$, that satisfies  \eqref{P6} for some $z\in\mathbb{R}$. There exists $\varepsilon_{0}>0$ only depending on the speed $c$, such that if $0<\varepsilon<\varepsilon_{0}$, then it holds
\begin{equation}
M^{3}-\frac{1}{4}E(u)M+\frac{1}{72}F(u)\le 0.
\label{P33}
\end{equation}
\end{Lem}\noindent
\textbf{Proof.}
The key is to show that $h\le 18 M$ on $\mathbb{R}$. Note that by \eqref{PP3} we know that 
$18M\ge c/4$ and that Lemma \ref{Lemma PP2} ensures that  $ \xi \in \Theta_z $ for $\varepsilon_{0}$ small enough. Let us set $\lambda=z-6.7$, $\mu=z+6.7$, and rewrite the function $h$ as
\begin{equation*}
  h(x)=\left\{
    \begin{aligned}
     &-v_{xx}-6v_{x}+16v,~~x<\lambda,\\
     &u-6v_{x}+12v,~~\lambda<x<\xi,\\
     &u+6v_{x}+12v,~~\xi<x<\mu,\\
     &-v_{xx}+6v_{x}+16v,~~x>\mu  \; . \\
    \end{aligned}
  \right.
  \end{equation*}
If $x\in\mathbb{R}\setminus\Theta_z$, using that $v_{xx}=4v-u$, \eqref{P14} and \eqref{P15}, it holds 
\begin{align*}
h&\le|v_{xx}|+6|v_{x}|+16v\le u+32v\le\frac{c}{9}\le 18M.
\end{align*}
If $\lambda<x<\xi$, then $v_{x}\ge 0$, and using that $u\le 6v$ on $\mathbb{R}$, we have 
$$
h=u-6v_{x}+12v\le 18v.
$$
If $\xi<x<\mu$, then $v_{x}\le 0$, and similarly using that $u\le 6v$ on $\mathbb{R}$, we get 
$$
h=u+6v_{x}+12v\le 18v.
$$
Therefore, it holds 
\begin{equation}
h(x)\le 18\max_{x\in\mathbb{R}}v(x)=18M,~~\forall x\in\mathbb{R}.
\label{P34}
\end{equation}

Now, combining \eqref{GG2}, \eqref{HH2} and \eqref{P34}, we get 
$$
F(u)-144M^{3}=\int_{\mathbb{R}}h(x)g^{2}(x)dx
\le\|h\|_{L^{\infty}(\mathbb{R})}\int_{\mathbb{R}}g^{2}(x)dx
\le 18M(E(u)-12M^{2}),
$$
and we obtain the lemma.\hfill $ \square$ \vspace{2mm}

\noindent
\textbf{Proof of Theorem \ref{Theorem P}.}
We argue as El Dika and Molinet in \cite{MR2542735}.
As noticed after the statement of the theorem, it suffices to prove \eqref{P7} assuming that $u \in L^{2}(\mathbb{R}) $ satisfies \eqref{P1}, \eqref{P2} and \eqref{P6}.  We recall that  $M=v(\xi)=\max_{x\in\mathbb{R}}v(x)$ and we set $\delta=c/6-M$. We first remark that if $\delta\le 0$, combining \eqref{P4} and \eqref{P8}, it holds
\begin{equation*}
\|u-\varphi_{c}(\cdot-\xi)\|_{\mathcal{H}}\le|E(u_{0})-E(\varphi_{c})|^{1/2}
\le O(\varepsilon),
\end{equation*}
that yields the desired result. 
Now suppose that $\delta>0$, that is the maximum of the function $v$ is less than the maximum of $\rho_{c}$. 
 Combining \eqref{P4}, \eqref{P5} and \eqref{P33}, we get
$$M^{3}-\frac{1}{4}E(\varphi_{c})M+\frac{1}{72}F(\varphi_{c})\le O(\varepsilon^{2}).$$
Using that $E(\varphi_{c})=c^{2}/3$ and $F(\varphi_{c})=2c^{3}/3$, our inequality becomes
$$\left(M-\frac{c}{6}\right)^{2}\left(M+\frac{c}{3}\right)\le  O(\varepsilon^{2}).$$
Next, substituting $M$ by $c/6-\delta$ and using that $[M+c/3]^{-1}<3/c$, we obtain
\begin{equation}
\delta^{2}\le O(\varepsilon^{2})\Rightarrow\delta\le O(\varepsilon).
\label{P35}
\end{equation}
Finally, combining \eqref{P4}, \eqref{P8} and \eqref{P35}, we infer that 
$$\|u-\varphi_{c}(\cdot-\xi)\|_{\mathcal{H}}\le C\sqrt{\varepsilon},$$
where $C>0$ only depends on the speed $c$. This completes the proof of  the stability of peakons.

\begin{Ack}
\normalfont
The author would like to thank his PhD advisor Luc Molinet for his help and his careful reading of this manuscript.
\end{Ack}

\bibliographystyle{plain}
\bibliography{mabiblio}

\end{document}